\magnification=\magstep1
\input amstex
\UseAMSsymbols
\input pictex
\vsize=23truecm
\NoBlackBoxes
\pageno=1

   \font\rmk=cmr8      \font\ttk=cmtt8


     \def\arr#1#2{\arrow <1.5mm> [0.25,0.75] from #1 to #2}

     

	       \centerline{\bf The representations of quivers of type $\Bbb A_n$. A fast approach.}
	       \bigskip
\centerline{Claus Michael Ringel}
		  \bigskip
It is well-known that a quiver $Q$ of type $\Bbb A_n$ is representation-finite,
and that its indecomposable representations are thin (all Jordan-H\"older multiplicities are 0 or 1). 
By now, various methods of proof are
known. The aim of this note is to provide a straight-forward arrangement of possible
arguments in order to avoid indices and clumsy inductive considerations, but also avoiding
somewhat fancy tools such as the Bernstein-Gelfand-Ponomarev reflection functors or
bilinear forms and root systems. The proof we present deals with representations in
general, not only finite-dimensional ones. 
We only will use first year linear algebra, namely
the existence of bases of vector spaces $V,W$ compatible with
a given linear map $V\to W$, and the existence of a basis of a vector space 
which is compatible with two given subspaces. 
      \medskip 
{\bf Theorem.} {\it Any representation of a quiver of type $\Bbb A_n$ is a direct sum of
thin representations.}
     \medskip 
Proof. If $n=2$, then we deal  with a linear map $f\:V \to W.$ Any first year
linear algebra course shows how to obtain a direct decomposition: Take a basis $\Cal B$
of the kernel of $f$, extend it by a family $\Cal B'$ to a basis of $V$. Now
$\{f(b')\mid b'\in \Cal B'\}$ is a basis of the image $f(V)$ and we extend it 
by a family $\Cal B''$ to a basis of $W.$
  
Thus, let $n\ge 3$ and  use induction. We deal with a quiver $Q$ with underlying graph
$$
{\beginpicture
\setcoordinatesystem units <1cm,1cm>
\multiput{$\circ$} at 0 0  1 0  2 0  4 0  5 0 /
\plot 0.2 0  0.8 0 /
\plot 1.2 0  1.8 0 /
\plot 2.2 0  2.4 0 /
\plot 3.6 0  3.8 0 /
\plot 4.2 0  4.8 0 /
\put{$\dots$} at 3 0
\put{$\ssize 1\strut$} at 0 -0.2
\put{$\ssize 2\strut$} at 1 -0.2
\put{$\ssize 3\strut$} at 2 -0.2
\put{$\ssize n-1\strut$} at 4 -0.2
\put{$\ssize n\strut$} at 5 -0.2

\endpicture}
$$ 
Let $M$ be a representation of $Q$ and $x$ a vertex of $Q$.
We call $x$ a {\it peak for $M$} provided for any arrow $\alpha\:y \to z$ 
the map $M_\alpha$ is injective in case $d(x,z) = d(x,y)-1$, and surjective in
case $d(x,z) = d(x,y)+1.$ Obviously, a thin indecomposable representation $M$ of $Q$
with $M_x \neq 0$ has $x$ as a peak. Second, a direct sum of modules with $x$ as
a peak, has $x$ as a peak. And third, if $x$ is a peak for $M$, then also 
for any direct summand of $M$.
    
We assume now by induction that all representations of quivers 
of type $\Bbb A_{n-1}$ are direct sums of thin representations. We first show: 
   \medskip 
(1) {\it Given a vertex $1 < x < n$ of $Q$, then
any representation $A$ of $Q$ can be written as a direct sum $A = B\oplus C\oplus D$,
where $B$ has $x$ as a peak, the support of $C$ is contained in $\{1,2,\dots,x-1\}$
and the support of $D$ is contained in $\{x+1,x+2,\dots,n\}$.}
    \medskip
Proof: We first look at the restriction $A'$ of $A$ to the subquiver $Q'$ with vertices $1,2,\dots,x$.
By assumption, we write $A'$ as a direct sum of thin indecomposable modules, 
say $A' = B'\oplus C'$,
where $B'$ is a direct sum of thin indecomposable representations of $Q'$ with
$x$ a peak
and $C'$ a direct sum of thin indecomposable representations of $Q'$ with $C'_x = 0.$
Second, we look at the restriction $A''$ of $A$ to the subquiver $Q''$ 
with vertices $x,x+1,\dots,n$.
Again by assumption, we write $A''$ as a direct sum of thin indecomposable modules, 
say $A'' = B''\oplus D''$,
where $B''$ is a direct sum of thin indecomposable representations of $Q''$ with $x$ a peak
and $D''$ a direct sum of thin indecomposable representations of $Q''$ with $D''_x = 0.$ 
Since $C'_x = 0 = D''_x$, we see that $A_x = B'_x = B''_x.$
Let $B$ be the subrepresentation of $A$ defined as follows: Its restriction to $Q'$ is
$B'$, its restriction of $Q''$ is $B''$.
Let $C$ be the 
subrepresentation of $A$ whose restriction to $Q'$ is
$C'$ and $C_y = 0$ for $y \ge  x$. 
Let $D$ be the 
subrepresentation of $A$ whose restriction to $Q''$ is
$D''$ and $D_y = 0$ for $y \le  x$. 
Then clearly $A = B\oplus C\oplus D$, with $B$ having $x$ as a peak, 
the support of $C$ is contained in $\{1,2,\dots,x-1\}$
and the support of $D$ is contained in $\{x+1,x+2,\dots,n\}$.
    \medskip
Using (1) twice, first for $x = 2$, then for $x = n\!-\!1$, we see:
      \medskip
(2) {\it Any representation $A$ of $Q$ can be written as a direct sum $A = B\oplus C \oplus D \oplus E$,
where $B$ has both $2$ and $n\!-\!1$ as peaks, the support of $C$ is contained
in $\{1\}$, the support of $D$ is contained in $\{3,4,\dots,n-2\}$ and the support of
$E$ is contained in $\{n\}.$}
    \medskip
Since by induction the representations $C,D,E$ are direct sums of thin representations,
we may assume that $A = B$, thus that $A$ has both $2$ and $n\!-\!1$ as peaks.
But if $2$ and $n\!-\!1$ both are peaks, all the maps $A_\alpha$ with $\alpha$ an
arrow in-between $2$ and $n\!-\!1$ are isomorphisms, thus, $A$ is isomorphic to a
representation where all these maps $A_\alpha$ are identity maps. 
It remains 
to look at the case $n=3$ and a representation with peak $2$. We have to show:
   \medskip
{\it Let $Q$ be a quiver of type $\Bbb A_3$ with graph $1 - 2 - 3$. Any representation of
$Q$ with peak $2$ is the direct sum of thin representations.}
Three different
orientations have to be discussed:
$$
{\beginpicture
\setcoordinatesystem units <1cm,1cm>
\put{\beginpicture
\multiput{$\circ$} at 0 0  1 0  2 0 /
\arr{0.2 0}{0.8 0}
\arr{1.8 0}{1.2 0}
\put{$\alpha$} at 0.5 0.2
\put{$\beta$} at 1.5 0.2
\put{$\ssize 1$} at 0 -0.25
\put{$\ssize 2$} at 1 -0.25
\put{$\ssize 3$} at 2 -0.25
\endpicture} at 0 0
\put{\beginpicture
\multiput{$\circ$} at 0 0  1 0  2 0 /
\arr{0.2 0}{0.8 0}
\arr{1.2 0}{1.8 0}
\put{$\alpha$} at 0.5 0.2
\put{$\beta$} at 1.5 0.2
\put{$\ssize 1$} at 0 -0.25
\put{$\ssize 2$} at 1 -0.25
\put{$\ssize 3$} at 2 -0.25
\endpicture} at 4 0
\put{\beginpicture
\multiput{$\circ$} at 0 0  1 0  2 0 /
\arr{0.8 0}{0.2 0}
\arr{1.2 0}{1.8 0}
\put{$\alpha$} at 0.5 0.2
\put{$\beta$} at 1.5 0.2
\put{$\ssize 1$} at 0 -0.25
\put{$\ssize 2$} at 1 -0.25
\put{$\ssize 3$} at 2 -0.25
\endpicture} at 8 0
\endpicture}
$$
In the first case, we deal with a representation $A$ such that both $A_\alpha$
and $A_\beta$ are injective, thus up to isomorphism we may assume that 
$A_\alpha$ and $A_\beta$ are inclusions of subspaces. Again 
any first year linear algebra course shows how to obtain a direct decomposition: 
Take a basis $\Cal B$ of the intersection $A_1\cap A_3$, extend it by a family
$\Cal B'$ to a basis of $A_1$ and by a family $\Cal B''$ to a basis of $A_3$.
Then the disjoint union 
$\Cal B\cup\Cal B'\cup \Cal B''$ is a basis of $A_1+A_2$ and we can extend this
by a family $\Cal B'''$ to obtain a basis of $A_2.$

In the second case, we deal with a representation $A$ such that $A_\alpha$
is injective, $A_\beta$ is surjective. Thus up to isomorphism we may assume that 
$A_\alpha$ is the inclusion of a subspace and we denote by $A'_3$ the kernel of 
$A_\beta$.  
Similar to the first case, we construct a basis of $A_2$
which is compatible with the two subspaces $A_1$ and $A'_3$, this yields the
required direct decomposition of $A$.

Finally, in the third case, we deal with a representation $A$ such that both $A_\alpha$
and $A_\beta$ are surjective. Let $A'_1$ be the kernel of $A_\alpha$ and
$A'_3$ the kernel of $A_\beta$. Again, as before,
we construct a basis of $A_2$
which is compatible with the two subspaces $A'_1$ and $A'_3$.
This completes the proof. 
     \bigskip
Of course, as a consequence we obtain: {\it Given two filtrations
$U_1\subseteq U_2 \subseteq \cdots \subseteq U_m$ and 
$U'_1\subseteq U'_2 \subseteq \cdots \subseteq U'_{m'}$ of a vector space $V$,
there is a basis of $V$ which is compatible with all these subspaces $U_i, U'_j.$}
Namely, look at the corresponding representation $A$
of the following quiver of type $\Bbb A_{m+m'+1}$
$$
{\beginpicture
\setcoordinatesystem units <1cm,1cm>
\multiput{$\circ$} at 0 0  1 0  3 0  4 0  5 0  7 0  8 0 /
\arr{0.2 0}{0.8 0}
\arr{3.2 0}{3.8 0}
\plot 1.2 0  1.6 0 /
\put{$\cdots$} at 2 0 
\arr{2.4 0}{2.8 0}

\arr{7.8 0}{7.2 0}
\arr{4.8 0}{4.2 0}
\plot 6.8 0  6.4 0 /
\put{$\cdots$} at 6 0 
\arr{5.6 0}{5.2 0}
\put{$\ssize 1$} at 0 -.3
\put{$\ssize 2$} at 1 -.3
\put{$\ssize m$} at 3 -.3
\put{$\ssize \omega$} at 4 -.3
\put{$\ssize m'$} at 5 -.3
\put{$\ssize 2'$} at 7 -.3
\put{$\ssize 1'$} at 8 -.3
\endpicture}
$$
 with $A_\omega = V$, $A_i = U_i$ for $1\le i \le m$ and $A_{j'} = U'_j$ for $1\le j\le m'$.
 \bigskip\bigskip
{\rmk 
e-mail: \ttk ringel\@math.uni-bielefeld.de \par}

\bye